\documentclass{article}
\usepackage{ijcai23}

\usepackage{times}
\usepackage{soul}
\usepackage{url}
\usepackage[hidelinks]{hyperref}
\usepackage[utf8]{inputenc}
\usepackage[small]{caption}
\usepackage{graphicx}
\usepackage{amsmath,amsfonts,amssymb, mathtools}
\usepackage{amsthm}

\usepackage{booktabs}
\usepackage{algorithm}
\usepackage{algorithmic}
\usepackage{makecell}
\usepackage[switch]{lineno}

\usepackage[table]{xcolor}    
\usepackage{mathtools}

\DeclarePairedDelimiter\floor{\lfloor}{\rfloor}

\title{{Machine Learning for Cutting Planes in Integer Programming: A Survey}}

\author{
Arnaud Deza$^1$
\and
Elias B. Khalil$^{1,2}$
\affiliations
$^1$Department of Mechanical \& Industrial Engineering, University of Toronto\\
$^2$SCALE AI Research in Data-Driven Algorithms for Modern Supply Chains
\emails
arnaud.deza@mail.utoronto.ca,
khalil@mie.utoronto.ca
}

\begin{document}
\maketitle
\begin{abstract}
    We survey recent work on machine learning (ML) techniques for selecting cutting planes (or \textit{cuts}) in mixed-integer linear programming (MILP). Despite the availability of various classes of cuts, the task of choosing a set of cuts to add to the linear programming (LP) relaxation at a given node of the branch-and-bound (B\&B) tree has defied both formal and heuristic solutions to date. ML offers a promising approach for improving the cut selection process by using data to identify promising cuts that accelerate the solution of MILP instances. This paper presents an overview of the topic, highlighting recent advances in the literature, common approaches to data collection, evaluation, and ML model architectures. We analyze the empirical results in the literature in an attempt to quantify the progress that has been made and conclude by suggesting avenues for future research.
\end{abstract}

\section{Introduction}
ML has recently been applied to accelerate the solution of optimization problems, with MILP being one of the most active research areas~\cite{bengio,kotary2021end,mazyavkina2021reinforcement}. A MILP is an optimization problem that involves both continuous and discrete variables, and aims to minimize or maximize a linear objective function $\boldsymbol{c}^{\intercal}\boldsymbol{x}$, over its decision variables $\boldsymbol{x} \in \mathbb{Z}^{|\mathbb{J}|} \times \mathbb{R}^{n - |\mathbb{J}|}$ while satisfying a set of $m$ linear constraints $ \boldsymbol{A}\boldsymbol{x}\leq \boldsymbol{b}$. Here, $\mathbb{J} \subseteq \{1,\cdots,n\}, |\mathbb{J}|\geq 1$, corresponds to the set of indices of integer variables. Similarly, Integer programming (IP) problems are of the same form only with discrete variables, i.e $\boldsymbol{x} \in \mathbb{Z}^{n}$. The MILP problem is written as:
\begin{align}
    \label{eqn:MILP}
    z^{IP} = \text{min}\{\boldsymbol{c}^{\intercal}\boldsymbol{x} \ \ | \ \   \boldsymbol{A}\boldsymbol{x}\leq \boldsymbol{b},  \ \    \boldsymbol{x} \in \mathbb{Z}^{|\mathbb{J}|} \times \mathbb{R}^{n - |\mathbb{J}|}  \}
\end{align}
The MILP formalism is widely used in supply chain and logistics, production planning, etc. While the MILP~\eqref{eqn:MILP} problem is NP-hard in general, modern solvers are able to effectively tackle large-scale instances, often to global optimality, using a combination of exact search and heuristic techniques. The backbone of MILP solving is the implementation of a tree search algorithm, \emph{Branch and Bound} (B\&B) \cite{b_and_b}, which relies on repeatedly solving computationally tractable versions of the original problem where discrete variables are relaxed to take on continuous values. Formally, we denote the linear programming (LP) relaxation of problem~\eqref{eqn:MILP}:
\begin{align}
    \label{eqn:LP}
    z^{LP} = \text{min}\{\boldsymbol{c}^{\intercal}\boldsymbol{x} \ \ | \ \   \boldsymbol{A}\boldsymbol{x}\leq \boldsymbol{b},  \ \    \boldsymbol{x} \in  \mathbb{R}^{n}  \}
\end{align}

\begin{figure}[t]
    \centering
\resizebox{\columnwidth}{!}{
    \includegraphics[scale = 0.39]{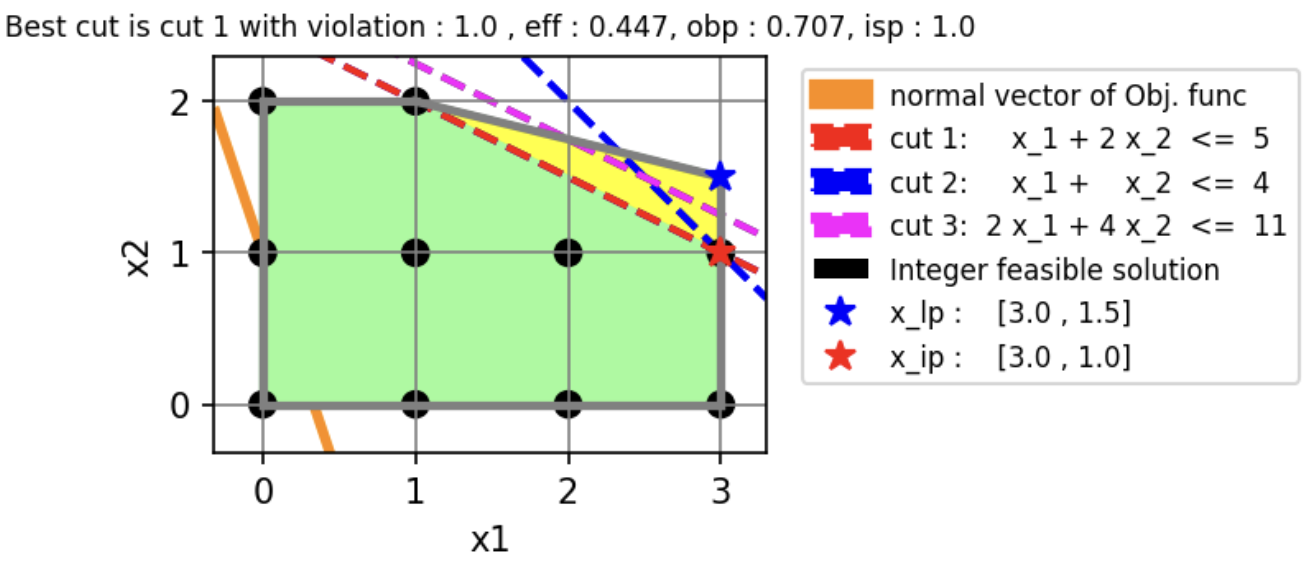}
}    \caption{A 2-dimensional IP. The optimum of the LP relaxation is shown as a blue star, whereas the integer optimum is shown as a red star. The colored cuts separate the LP optimum as desired. The best cut is cut 1 as it produces the convex hull, shaded in pale green, and evaluation metrics calculated for this cut are shown above the graph.}
    \label{fig:cut_2d}
\end{figure}

To render the B\&B search more efficient, valid linear inequalities (or tightening constraints) to problem~\eqref{eqn:MILP}-- \textit{cutting planes} --are added to LP relaxations of the MILP with the aim of producing tighter relaxations and thus better lower bounds to problem~\eqref{eqn:MILP}, as illustrated in figure \ref{fig:cut_2d}; this approach is referred to as the Branch and Cut (B\&C) algorithm. Cuts are essential for MILP solving, as they can significantly reduce the feasible region of the B\&B algorithm and exploit the structure of the underlying combinatorial problem, which a pure B\&B approach does not do. The incorporation of cutting planes in the B\&B search necessitates appropriate filtering and selection as there can be many such valid inequalities and adding them to a MILP comes at a cost in computation time. As such, \textit{cut selection} has been an area of active research in recent years.

Various families of general-purpose and problem-specific cuts have been studied theoretically and implemented in modern-day solvers~\cite{Santanu}. However, there is an overall lack of scientific understanding regarding many of the key design decisions when it comes to incorporating cuts in B\&B. Traditionally, the management of cutting plane generation and selection is governed by hard-coded parameters and manually-designed expert heuristic rules that are based on limited empirical results. These rules include deciding:
\begin{itemize}
    \item[--] the number of cuts to generate and select;
    \item[--] the number of cut generation and selection rounds;
    \item[--] add cuts at the root node or all nodes of the B\&B tree;
    \item[--] which metrics to use to score cuts for selection;
    \item[--] which cuts to remove and when to remove them.
\end{itemize}
A cut selection strategy common in modern solvers is to use a scoring function parameterized by a linear weighted sum of cut quality metrics that aim to gauge the potential effectiveness of a single cut. This process is iteratively used to produce a ranking among a set of candidate cuts. However, manually-designed heuristics and fixed weights may not be optimal for all MILP problems \cite{ACS}, and as such, researchers have proposed using ML to aid in cuts selection. Figure \ref{fig:cut_2d} gives the reader a 2D visualization of cuts and potential selection rules.

Recently, the field of ML for cutting plane selection has gained significant attention in MILP~\cite{Columbia,huawei,tree_complexity,local_cuts,ACS,look_ahead,Analytic_Centers} and even mixed-integer \textit{quadratic} programming~\cite{BalteanLugojan2019ScoringPS} and \textit{stochastic} integer programming~\cite{Benders_cut}. To select more effective cutting planes, researchers have proposed ML methods ranging from reinforcement to imitation and supervised learning. This survey aims to provide an overview of the field of ML for MILP cut selection starting with the relevant ML and MILP background (Section~\ref{sec:background}), the state of cut selection in current MILP solvers (Section~\ref{sec:cutsel}), recently proposed ML approaches (Section~\ref{sec:ml}), relevant learning theory (Section~\ref{sec:theory}), and avenues for future research (Section~\ref{sec:conclusion}). 

\section{Background}
\label{sec:background}

\subsection{Integer programming and Cutting planes}
Cutting planes (or \emph{cuts}) are valid linear inequalities to problem~\eqref{eqn:MILP} of the form $\mathbf{\boldsymbol{\alpha}}^T x \leq \beta, \mathbf{\boldsymbol{\alpha}} \in \mathbb{R}^n, \beta \in \mathbb{R}$. They are ``valid" in the sense that adding them to the constraints of the LP relaxation in ~\eqref{eqn:LP} is guaranteed not to cut off any feasible solutions to~\eqref{eqn:MILP}. Additionally, one seeks cuts that separate the current LP solution, $x^*_{LP}$, from the convex hull of integer-feasible solutions; see Figure~\ref{fig:cut_2d}. While adding more cuts can help achieve tighter relaxations in principle, a clear trade-off exists: as more cuts are added, the size of the LP relaxation grows resulting in an increased cost in LP solving at the nodes of the B\&B tree~\cite{Tobias_thesis}. Adding too few cuts, however, may lead to a large number of nodes in the search tree as more  branching is required.

We note that the so-called \textit{cutting plane method} can theoretically  solve integer linear programs by iteratively solving relaxed versions of the given problem then adding cuts to separate the fractional relaxed solution $x^*_{LP}\in \mathbb{R}^n$, and terminating when $x^*_{LP}\in \mathbb{Z}^n$. Despite theoretical finite convergence results for the cutting plane method using Gomory cuts, numerical errors and design decisions such as cut selection will often prevent convergence to an optimal solution in practice.

\subsection{ML for Combinatorial Optimization }
The use of ML in MILP and combinatorial optimization (CO) has recently seen some success, with a diversity of approaches in the literature \cite{bengio} that can be split into two main categories. The first is to directly predict near-optimal solutions conditioned on the representation of particular instances. Notable examples of this include learning for quadratic assignment problems\cite{nowak}, solving CO problems with pointer networks \cite{Vinyals}, and using attention networks to solve travelling salesman problems \cite{att_tsp}. These approaches aim to completely replace traditional solvers with ML models and are hence very appealing given their black-box nature. In contrast, a second approach focuses on automating decision-making in solvers through learned inductive biases. This can take on the form of replacing certain challenging algorithmic computations with rapid ML approximations or using newly learnt heuristics to optimize solution time. Notable examples include the learning of computationally challenging variable selection rules for B\&B \cite{khalil2016learning,GCNN,zarpellon2021parameterizing}, learning to schedule heuristics \cite{learn_to_schedule}, or learning variables biases \cite{mip_gnn}.

\noindent\textbf{Representing MILPs for ML.} One of the key challenges in applying ML to MILP is the need for efficient and clever feature engineering. This requires a deep understanding of solver details, as well as a thorough understanding of the underlying problem structure. In recent years, graph neural networks (GNNs) have emerged as a popular architecture for several ML applications for MILP~\cite{cappart}. GNNs have the ability to handle sparse MILP instances and exhibit permutation invariance, making them well-suited for representing MILP instances. The GNN operates on the so-called \textit{variable-constraint graph} (VCG) of a MILP. The VCG has $n$ variable nodes and $m$ constraint nodes corresponding to the decision variables and constraints of~\ref{eqn:MILP}. The edges between a variable node $j$ and constraint node $k$ represent the presence of that variable $x_j$ in constraint $k$ (i.e, $A_{jk} \neq 0$), where the weight of the edge is $A_{jk}$. 


\subsection{Common Learning Paradigms in ML for CO}

\noindent\emph{\textbf{Supervised Learning}}: The simplest and most common learning paradigm is supervised learning (SL) where the learning algorithm aims to find a function (model) $f :X \rightarrow Y ,\  f \in F$, where $F$ is the function's hypothesis space, given a labelled dataset of $N$ training samples of the form $\{(x_1, y_1),\cdots(x_N, y_N)\}$ where $x_i\in X$ represents the feature vector of the i-th sample and $y_i\in Y$ its label. The goal is to find a $f$ such that a loss function $L(y_i,\hat{y_i})$ measuring how well predictions, $\hat{y_i}$, from $f$ \emph{fit} to the training data with the hope of generalizing to unseen test instances.

\noindent\emph{\textbf{Reinforcement Learning}}: A Markov Decision Process (MDP) is a mathematical framework for modelling sequential decision-making problems commonly used in reinforcement learning (RL). At time step $t\geq0$ of an MDP, an agent in state $s_t \in \mathcal{S}$ takes an action $a_t\in \mathcal{A}$ and transitions to the next state $s_{t+1} \sim p(\cdot|s_t,a_t)$, incurring a scalar reward $r_t\in\mathbb{R}$. A policy, denoted by $\pi$, provides a mapping $\mathcal{S} \mapsto P(\mathcal{A})$ from any state to a distribution over actions $\pi(\cdot|s_t)$. The goal of an MDP is to find a policy $\pi$ that maximises the expected cumulative rewards over a horizon $T$, i.e, $max_{\pi}$ J($\pi$) = $\mathbb{E}_{\pi}[\sum_{t=0}^{T-1}\gamma^tr_t]$,  where $\gamma \in (0,1]$ is the discount factor. 

\noindent\emph{\textbf{Imitation Learning}}: Imitation Learning (IL) aims to find a policy $\pi$ that mimics the behaviour of an expert in a given task through demonstrations. This is often formalized as an optimization problem of the form $\min_{\pi} L(\pi,\hat{\pi})$, where $\hat{\pi}$ is the expert policy and $L$ is a loss function measuring the difference between the expert and the learned policy. IL can be seen as a special case of RL where the agent's objective is to learn a policy that maximizes the likelihood of the expert's actions instead of the expected cumulative reward. In  ML-for-MILP, IL (and SL) has been used to amortize the cost of powerful yet computationally intractable scoring functions. Such functions appear in cut generation/selection~\cite{Amaldi2014CoordinatedCP,Coniglio} and have been recently approximated using IL~\cite{look_ahead}, as we will discuss hereafter.

\section{Cutting Planes in MILP solvers}
\label{sec:cutsel}

\subsection{Cut generation}
At a node of the search tree during the B\&C process and prior to branching, the solver runs the cutting plane method for a pre-specified \emph{number of separation rounds}, where each round $k$ involves i) solving a continuous relaxation $P^{(k)}$ to get a fractional $x^k_{LP}$; ii) generating cuts $\mathcal{C}^k$, referred to as the \emph{cut-pool} followed by selecting a subset $\mathcal{S}^{k}\subseteq\mathcal{C}^k$; iii) adding $S^{k}$ to the relaxation and proceeding to round $k+1$. After $k$ separation rounds the LP relaxation $P^k$ will consist of the original constraints $\boldsymbol{A}\boldsymbol{x}\leq \boldsymbol{b}$ and any cuts of the form ($\mathbf{\boldsymbol{\alpha}},\beta$) that have been selected.
Concretely, we write $P^k$ as:
\begin{equation}
\begin{aligned}
    P^k &=\{ \boldsymbol{A}\boldsymbol{x}\leq \boldsymbol{b},  \ \mathbf{\boldsymbol{\alpha}}^T \boldsymbol{x} \leq \beta \ \forall \ (\mathbf{\boldsymbol{\alpha}},\beta) \in \bigcup_{i=1}^{k} S^i,  \ \    \boldsymbol{x} \in  \mathbb{R}^{n}  \}\\
     &=\{\boldsymbol{A}^k\boldsymbol{x}\leq \boldsymbol{b}^k, \boldsymbol{x} \in  \mathbb{R}^{n}  \} \quad \text{with  } \ \boldsymbol{A}^0=\boldsymbol{A}, \  \boldsymbol{b}^0=\boldsymbol{b}
\end{aligned}
\end{equation}

\noindent \emph{Global cuts} are cuts that are generated at the root node whereas \emph{local cuts} are cuts generated at nodes further down the tree. Traditionally, solely relying on global cuts is referred to as \emph{Cut \& Branch}, in comparison to B\&C which uses both global and local cuts. Solvers use various hard-coded parameters determined experimentally to control the number of separation rounds, types of generated cuts, their frequency, priority among separators, whether to use local cuts or not, among others. We use SCIP's internal cut selection subroutine to highlight key decisions regarding cut selection. However, similar methodologies are used in other MILP solvers.

\subsection{Cut selection}
During the cut selection phase, the solver selects a subset of cuts $S^k$ from the cut-pool to add to the MILP. In particular, SCIP scores each cut using a simple linear weighted sum of cut quality metrics from the MILP literature: 
\begin{equation}
    \label{eqn:ScipScore}
    S = \lambda_1 \textbf{eff} + \lambda_2  \textbf{dcd} + \lambda_3 \textbf{isp} + \lambda_4 \textbf{obp} \quad 
    \mathbf{\lambda} \geq  \mathbf{0}, \mathbf{\lambda} \in \mathbb{R}^4
\end{equation}
Here, the weights ($\lambda_1,\lambda_2,\lambda_3,\lambda_4$) are solver parameters that correspond to four metrics: efficacy (\textbf{eff}), directed cutoff distance (\textbf{dcd}), integer support (\textbf{isp}), and objective parallelism (\textbf{obp}). Cuts are then added greedily by the highest-ranking score $S$ and added to $S_k$, followed by filtering the remaining cuts for parallelism. This is done until a pre-specified number of cuts have been selected or no more candidate cuts remain. The current default weights for SCIP version 8.0 \cite{scip8} are $\mathbb{\lambda_{DEF}}^T = [0.9, 0.0, 0.1, 0.1]$.

\subsection{Cut Metrics}
Cheap metrics such as the ones above are used to gauge the potential bound improvement that will result from adding a cut  ($\boldsymbol{\alpha},\beta$) to a current relaxation with optimal solution $x_{LP}$ and best known incumbent $\hat{x}$. 
\emph{Efficacy} is the {Euclidean distance} between the hyperplane $\boldsymbol{\alpha}^T x = \beta$ and $x_{LP}$ and can be interpreted as the \emph{distance cut off by a cut} \cite{Wesselmann}, eXpressed as:
\begin{equation}
    \label{eqn:eff}
    \text{\textbf{eff}}(\boldsymbol{\alpha},\beta,x_{LP})\coloneqq \frac{\boldsymbol{\alpha}^T x_{LP} - \beta}{\|\boldsymbol{\alpha}\|}.
\end{equation}
\emph{Directed cutoff distance} \cite{scip_6} is the distance between the hyperplane $\boldsymbol{\alpha}^T x = \beta$ and $x_{LP}$ in the direction of $\hat{x}$ and is measured as:
\begin{equation}
    \label{eqn:dcd}
    \text{\textbf{dcd}}(\boldsymbol{\alpha},\beta,x_{LP},\hat{x})\coloneqq \frac{\boldsymbol{\alpha}^T x_{LP} - \beta}{|\boldsymbol{\alpha}^T\mathbf{y}|}, \ \mathbf{y}\coloneqq\frac{\hat{x}-x_{LP}}{\|\hat{x}-x_{LP}\|}.
\end{equation}
The \emph{support of a cut} is the fraction of coefficients $\boldsymbol{\alpha}_i$ that are non-zero; sparser cuts are preferred for computational efficiency and numerical stability~\cite{Santanu}. The integer support takes this notion one step further by considering the fraction of coefficients corresponding to integer variables that are non-zero, measured as:
\begin{equation}
    \label{eqn:isp}
    \text{\textbf{isp}}(\boldsymbol{\alpha}) \coloneqq \frac{\sum_{i\in\mathbb{J}}^{}\text{NZ}(\boldsymbol{\alpha}_i)}{\sum_{i=1}^{n}\text{NZ}(\boldsymbol{\alpha}_i)}, \ \text{NZ}(\boldsymbol{\alpha}_i) \coloneqq \begin{cases} 
          0 & \text{if } \boldsymbol{\alpha}_i=0 \\
          1 & \text{else}
       \end{cases} \\
\end{equation}
\emph{Objective parallelism} is measured by considering the cosine of the angle between $c$ and $\boldsymbol{\alpha}$, with $\textbf{obp}(\boldsymbol{\alpha},c) = 1$ for cuts parallel to the objective function:
\begin{equation}
    \label{eqn:obp}
    \text{\textbf{obp}}(\boldsymbol{\alpha},c)\coloneqq \frac{|\boldsymbol{\alpha}^T c|}{\|\boldsymbol{\alpha}\|  \|c\| }
\end{equation}

\noindent More directly useful, but expensive, evaluation metrics can be measured by solving the relaxation obtained by adding the selected cuts and observing its objective value. Specifically, the \emph{integrality gap} (IG) after separation round $k$ is given by the bound difference $g^k \coloneqq z_{IP} -z^{k}\geq 0$ whereas the integrality gap closed (IGC) is measured as :
\begin{equation}
    IGC^{k}\coloneqq\frac{g^0-g^k}{g^0}=\frac{z^k-z^0}{z^{IP}-z^0} \in [0,1]
\end{equation}
and represents the factor by which the integrality gap is closed between the first relaxation $P^0$ and the relaxation $P^k$ obtained after $k$ separation rounds \cite{Columbia}.

\subsection{Precursors to Learning to Cut}
Since cut selection is not an exact science and no formal guideline exists, traditional methods to find good cut selection parameters for Eq.~\eqref{eqn:ScipScore} rely on performing parameter sweeps using appropriately designed grid searches. For instance, the first large-scale computational experiment regarding cut selection was presented in \cite{Tobias_thesis} and many more computational studies in SCIP have since been published \cite{scip7,scip8}. \cite{Tobias_thesis} is the basis for SCIP's scoring function and involves testing many configurations of cut metrics presented in \cite{Wesselmann} for \ref{eqn:ScipScore} and demonstrates a large decrease in overall solution time and nodes in the B\&B tree if the parameters are tuned properly.

Overall, the many hard-coded parameters that are used in MILP solvers can be tuned either by MILP experts  or through black-box algorithm configuration methods. These  include  grid search, but also more sophisticated methods such as sequential model-based optimization methods \cite{SMAC3} or black-box local search~\cite{xu2011hydra}.

\begin{table*}[t]
    \centering
  \rowcolors{2}{gray!15}{white}
    \begin{tabular}{p{1.8cm}p{3.1cm}p{1.5cm}p{1.3cm}p{3.9cm}p{1.2cm}p{1.4cm}}
        \toprule
                \textbf{Paper}  & \textbf{Learning task}  & \textbf{Solver}  & \textbf{ML paradigm} & \textbf{Instance type/source} & \textbf{Instance size} & \textbf{ML model}  \\
        \midrule
        \midrule
        \cite{Columbia}     & score single cut & Gurobi & RL & Synthetic & small & Attention \& LSTM \\
        \cite{look_ahead}       & score single cut  & SCIP & IL & Synthetic + NN verification & large & GNN \\
        \cite{BalteanLugojan2019ScoringPS}   & score single cut & CPLEX & SL & QP + QCQP & large & MLP      \\
        \cite{Benders_cut}   & classify single cut & Bender's for 2SP & SL & CFLP + CMND & medium & SVM      \\
        
        \cite{huawei}       & score bag of cuts & Proprietary solver & SL & Proprietary data & large & MLP  \\
        
        \cite{ACS}          & learning weights for cut scoring & SCIP  & RL & MIPLIB 2017 & large & GNN      \\
        
        
        \cite{local_cuts}   & learning when to cut & Xpress & SL & MIPLIB 2017 + Proprietary benchmark & large & Random Forest \\
        \cite{hierarchical} & learn to score cuts and order/number of cuts & SCIP & RL & MIPLIB 2017 + Synthetic + Proprietary benchmark & large & LSTM \& Pointer\\
        \bottomrule
    \end{tabular}
    \caption{Table summarizing and categorizing tasks tackled by research embedding ML for cut management in optimization solvers. The three instance size classifications, small, medium, and large correspond to instances with $n \times m$ in the range $[0,1000],[1000,5000],[5000,\infty]$ respectively. Additionally, QCQP refers to quadratically constrained QPs and Synthetic refers to the 4 IP instances presented in \protect\cite{Columbia} which are integer/binary packing, max cut and production planning problems}
    \label{tab:table_TASKS}
\end{table*}

\section{Learning to Cut}
\label{sec:ml}

The research on ``\emph{Learning to Cut}" can be categorized along three axes: the choice of the cut-related learning task, the ML paradigm used, and the optimization problem class of interest (MILP or others). We use these axes to organize the survey. Table \ref{tab:table_TASKS} provides a classification of the surveyed papers.

\subsection{Directly scoring individual cuts in MILP}
\subsubsection*{Scoring using Reinforcement Learning~\cite{Columbia}}
\citeauthor{Columbia},~\citeyear{Columbia}, were the first to motivate and experimentally validate the use of \emph{any} {learning} for cut selection in MILP solvers. The authors present an MDP formulation of the iterative cutting plane method (discussed in Section~\ref{sec:background}) for Gomory cuts from the LP tableau. A single cut is selected in every round via a neural network (NN) that predicts cut scores and produces a corresponding ranking. Given that GNNs were still in their infancy at the time, the authors instead used a combination of attention networks for order-agnostic cut selection and an LSTM network for IP size invariance. The authors used evolutionary strategies as their learning algorithm and considered the following baseline selection policies: maximum violation, maximum normalized violation, lexicographical rule, and random selection.

At iteration $t$ of the proposed MDP, the state $s_t \in \mathcal{S}$ is defined by $\{\mathcal{C}^{(t)},x^*_{LP}(t),P^{(t)}\}$ and the discrete action space $\mathcal{A}$ includes available actions given by $\mathcal{C}^{(t)}$, i.e., the Gomory cuts parameterized by $\mathbf{\alpha_i} \in \mathcal{R}^n,\beta_i\in\mathcal{R} \ \forall \ i \in \{0,\dots,|\mathcal{C}^{(t)}|\}$ that could be added to the relaxation $P^{(t)}$. The reward $r_t$ is the objective value gap between two consecutive LP solutions, i.e., $r_t\coloneqq\mathbf{c}^{T}[x^*_{LP}(t+1) - x^*_{LP}(t)] \geq 0$ which when combined with a discount factor $\gamma < 1$, encourages the agent to reduce the IG and reach optimality as fast as possible. Given a state $s_t = \{\mathcal{C}^{(t)},x^*_{LP}(t),P^{(t)}\}$ and an action $a_t$ (i.e, a chosen Gomory cut $\mathbf{\alpha_i}^Tx\leq\beta_i$), the new state $s_{t+1} = \{\mathcal{C}^{(t+1)},c,x^*_{LP}(t+1),P^{(t+1)}\}$ is determined by i) solving the new relaxation $P^{(t+1)} = P^{(t)} \cup\{\alpha_i^Tx\leq \beta_i\}$ using the simplex method to get $x^*_{LP}(t+1)$, ii) generating the new set of Gomory cuts $\mathcal{C}^{(t+1)}$ read from the simplex tableau. 

The RL approach significantly outperformed all metrics by effectively closing the IG with the fewest number of cuts for four sets of synthetically generated IP instances. They demonstrated generalization across the instance types as well as across instance sizes in two test-bed environments: 1) pure cutting plane method, 2) $B\&C$ using Gurobi Callbacks. 

However, limitations of this work include weak baselines, the restriction to Gomory cuts and a state encoding that does not scale well for large-scale instances as the input to the NN includes all constraints and available cuts. Additionally, the instance sizes considered were fairly small for research in MILP which may have been acceptable given that this was early work in this space. Although many recent papers outperform this approach, it is significant given that it is the first paper that appropriately defines an RL task for cut selection in the cutting plane method or B\&C.

\subsubsection*{Scoring using Imitation Learning~\cite{look_ahead}}
In a recent paper, \citeauthor{look_ahead} demonstrate the strength of a greedy selection rule that explicitly looks ahead to select the cut that yields the best bound improvement, but they note that this approach is too expensive to be deployed in practice. They propose the lookahead score, $s_{LA}$, that measures the increase in LP relaxation value obtained from adding a cut $C_j$ to an LP relaxation $P$. Formally, $C_j \in \mathcal{C}$ where $\mathcal{C}$ is a pool of candidate cuts, and $C_j$ is parameterized by $(\boldsymbol{\alpha},\beta)$. Let $z^j$ denote the optimal value of LP relaxation $P^j \coloneqq P \cup \{\boldsymbol{\alpha}^T x \leq \beta \}$, the new relaxation obtained by adding $C_j$ to $P$. The (non-negative) lookahead score then reads:
\begin{align}
    s_{LA}(C_j,P) \coloneqq z^j - z.
\end{align}
In response to this metric's computational intractability, the authors propose a NN architecture, NeuralCut, that is trained using imitation learning with $s_{LA}$ as its expert. The prohibitive cost of the expensive lookahead, which requires solving one LP per cut to obtain $z^j$, is thus alleviated by an approximation of the score. 
The authors collect a dataset of expert samples by running the cut selection process for 10 rounds and recording the cuts, LP solutions, and scores from the Lookahead expert creating samples specified by $\{\mathcal{C},P,\{s_{LA}(C_j,P)\}_{C_j\in\mathcal{C}}\}$. They use this expert data to learn a scoring $\hat{s}$ that mimics the lookahead expert by minimizing a soft binary entropy loss overall cuts,
\begin{gather}
    L(s) \coloneqq -\frac{1}{|\mathcal{C}|}\sum_{C \in \mathcal{C}}^{}q_C \log s_C + (1-q_C)\log (1-s_C),
\end{gather}
where $q_C \coloneqq \frac{s_{LA}(C)}{s_{LA}(C^*_{LA})}$ and $C^*_{LA} = \text{argmax}_{C \in \mathcal{C}} s_{LA}(C)$. To encode the cut selection decision that is described by the cut-pool and LP relaxation pair $(\mathcal{C},P)$, the authors use a tripartite graph whose nodes hold feature vectors for variables, constraints of P and any cuts added.\\
The 4 synthetic IP instance classes from \cite{Columbia} were used in this work. Only large instances were used given that small and medium-sized instances were observed to be too easy to solve. To evaluate their approach, the GNN is deployed for 30 consecutive separation rounds and adds a single cut per round in a pure cutting plane setting. The results show that NeuralCut exhibits great generalization capabilities, a close approximation of the lookahead policy, and outperforms the approach in \cite{Columbia} as well as many of the manual heuristics in \cite{Wesselmann} for 3 out of the 4 instances types; the packing instances tied for many methods and did not significantly benefit from a lookahead scorer or NeuralCut. To stress-test their approach, the authors employ NeuralCut at the root node in B\&C for a challenging dataset of NN verification problems \cite{Vinod_Nair} which are harder for SCIP to solve due to their larger size and notoriously weak formulations. They demonstrated clear benefits to the learned cut scorer.

A drawback of $s_{LA}$ is its limitation to scoring \textit{a single cut} due to the computational intractability of scoring a subset of cuts, of which there are combinatorially many. Additionally, although \cite{look_ahead} improves on \cite{Columbia}, both approaches have the inherent flaw of scoring each cut independently without taking into account the \emph{collaborative nature} of the selected cuts (i.e, complementing each other and uniquely collaborating in tightening relaxations). 

\subsubsection*{Scoring using Hierarchical Sequence Models~\cite{hierarchical}}

The authors of the most recent paper on cut selection, \cite{hierarchical}, motivate two new factors to consider for efficient cut selection: the number and order of cuts. They introduce a hierarchical sequence model (HEM), trained using RL, consisting of a two-level model: (1) a higher-level LSTM network that learns the number of cuts to be selected by outputting a ratio $k \in [0,1]$ of the cuts to be selected, (2) a lower-level pointer network formulating cut selection as a sequence-to-sequence learning problem that ultimately generates an ordered subset of cuts of size $\floor{N*k}$, where $N$ is the size of a pool of cuts. Experiments show that HEM, in comparison to rule-based and learned baselines from \cite{look_ahead,Columbia}, improves solver runtime by 16.4\% and primal dual integral by 33.48\% on synthetic MILP problems. Additionally, HEM was deployed on challenging benchmark instances from MIPLIB 2017~\cite{gleixner2021miplib} and large-scale production planning problems, showing great generalization capabilities. The authors also visualize cuts selected by HEM and \cite{look_ahead,Columbia}, demonstrating that HEM captures the order information and the interaction among cuts, leading to improved selection of complementary cuts. The two-level model of HEM is the only learned approach to consider the collaborative nature and interactions among cuts for efficient cut selection which was previously only explored in the MILP literature \cite{Coniglio}.

\subsection{Directly scoring individual cuts for Non-convex Quadratic Programming and Stochastic Programming}

The first work to incorporate any type of \emph{learning} for data-driven cut selection policies, even prior to~\cite{Columbia}, is that in~\cite{BalteanLugojan2019ScoringPS}. It similarly focuses on estimating lookahead scores for cuts. The lookahead criterion in their setting, non-convex quadratic programming (QP), involves solving a semidefinite program which is not viable in a B\&C framework. Although a different optimization setting, many ideas from MILP still apply and a similar approach of employing a NN estimator that predicts the objective improvement of a cut is used. However, a supervised regression task is considered which resulted in a trained multilayer perceptron (MLP) that exhibited speed-ups for evaluating cut selection measures approximately on the order of $2$x, $30$x and $180$x when compared to LAPACK's eigendecomposition method \cite{anderson1999lapack}, Mosek solver \cite{aps2019mosek} and SeDuMi solver \cite{polik2007sedumi} respectively.

Another optimization setting where appropriate cut selection is crucial is two-stage stochastic programming (2SP)~\cite{ahmed2010two}. Traditional solution techniques to 2SP include using Bender's decomposition which leverages problem structure through objective function approximations and the addition of cuts to sub-problem relaxations and a relaxed master problem. The authors in \cite{Benders_cut} leverage SL to train support vector machines (SVM) for the binary classification of the usefulness of a Bender's Cut and observe that their model allows for a reduction in the total solving time for a variety of 2SP instances. More specifically, solution time reductions ranging from $6\%$ to $47.5\%$ were observed on test instances of capacitated facility location problems (CFLP) and slightly smaller reductions were observed for multi-commodity network design (CMND) problems. 

\subsection{Directly scoring a bag of cuts for MILP}
In contrast to learning to score individual cuts, the authors in \cite{huawei} tackle cut selection through \textit{multiple instance learning} \cite{MIL} where they train a NN, in a supervised fashion, to score a \emph{bag} of cuts for Huawei's proprietary commercial solver. More specifically, the training samples, denoted by the tuple $\{\mathcal{P},C^\prime,r\}$, are collected using active and random sampling \cite{bello2016neural} where $r$ is the reduction ratio of solution time for a given MILP, with relaxation $\mathcal{P}$, when adding the bag of cuts $C^\prime$. The authors formulate the learning task as a binary classification problem, where the label of a bag $C^\prime$ is $1$ if $r$ ranks in the top $\rho\%$ highest reduction ratios for a given MILP, ($0$ otherwise), and $\rho \in (0,100)$ is a tunable hyper-parameter controlling the percentage of positive samples. At test time, the NN is used to assign scores to all candidate cuts and then select the top $K\%$ cuts with the highest predicted scores, where $K$ is another hyper-parameter. Other notable design decisions include designing bag features from aggregated cut features and a cross-entropy loss with L2 regularization to combat overfitting. 

The data for this work consisted of synthetic MILP problems solved within 25 seconds and large real-world production planning problems ranging in the $10^7$ variables. The results clearly demonstrate the benefit of a learned scorer by comparing their approach to rules from \cite{Wesselmann} and to a fine-tuned manual selection rule used by Huawei's proprietary solver. Once again, this method suffers from fixing the size/ratio of selected cuts, $K$, and scores cuts independently which neglects the preferred collaborative nature of selected cuts. Note that this is despite training the model to predict the quality of a bag of cuts: at test time, a ``bag" has only a single cut. 

\subsection{Learning Adaptive Cut Selection Parameters}
Rather than directly predicting cut scores,~\citeauthor{ACS},~\citeyear{ACS}, motivate learning instance-specific weights, $\mathbb{\lambda_{ACS}} \in \mathbb{R}^4$, for the SCIP cut scoring function in Eq.~\eqref{eqn:ScipScore}. 
The goal is to improve over the default parameterization $\mathbb{\lambda_{DEF}}$ in terms of relative gap improvement (RGI)  at the root node LP with 50 separations rounds of 10 cuts per round and a best-known primal bound.
Besides the learning approach proposed by the authors, a grid search over convex combinations of the four weights, $\sum_{i=1}^{4}\lambda_i = 1$, where $\lambda_i = \frac{\beta_i}{10}, \beta_i \in \mathbb{N}$, was performed individually for a large subset of MIPLIB 2017 instances. This experiment demonstrates the potential for improvement that one could get with instance-specific weights. The resulting parameters, referred to as $\mathbb{\lambda_{GS}}$, resulted in a median RGI of 9.6\%. 

The GNN architecture and VCG features are based on \cite{GCNN} but LP solution features are not used. The output of the model $\mu \in \mathcal{R}^4$ represents the mean of a multivariate normal distribution $\mathbb{N}_4(\mu,\gamma \mathbb{I})$, with $\gamma \in \mathbb{R}$ (a hyper-parameter) that is sampled to generate instance-specific parameters, $\mathbb{\lambda_{ACS}}$. Although the authors claim to use RL to train their GNNs, they fail to appropriately define the sequential nature of their MDP given that the time horizon, $T$, is 1. In their MDP, an action $a_t$ corresponds to a weight configuration $\mathbb{\lambda_{ACS}}$ sampled from $\mathbb{N}_4(\mu,\gamma \mathbb{I})$ which will in turn result in an RGI that is used as the instant reward $r_t$. As such, we consider their work to belong to instance-specific algorithm configuration~\cite{malitsky2014instance}, and the gradient descent approach used to train the GNN can be seen as approximating the unknowing mapping from (instance, parameter configuration) to RGI. 

GNN trained on an individual instance basis were able to achieve a relative gap improvement of 4.18\%. However, when trained over MIPLIB 2017 itself, a median RGI of 1.75\%  was achieved whereas a randomly initialized GNN produced an RGI of 0.5\%. The authors also observe that $\mathbb{\lambda_{GS}}$ and $\mathbb{\lambda_{ACS}}$  differ heavily from $\mathbb{\lambda_{DEF}}$ as seen in Table~\ref{tab:lambdas} and none of the values tend towards zero, meaning that all the metrics are able to provide utility depending on the given instance.
\begin{table}[h]
    \centering
    \resizebox{\columnwidth}{!}{%
    \begin{tabular}{c|ccc|ccc}
    \toprule
    Method & \multicolumn{3}{c}{Grid Search $\mathbb{\lambda_{GS}}$} & \multicolumn{3}{c}{ACS approach $\mathbb{\lambda_{ACS}}$} \\
    \midrule
    Parameter & Mean & Median & Std. Dev & Mean & Median & Std. Dev \\
    \midrule
    $\lambda_1$ (\textbf{eff}) & 0.179 & 0.100 & 0.216 & 0.294 & 0.286 & 0.122\\
    \midrule
    $\lambda_2$ (\textbf{dcd}) & 0.241 & 0.200 & 0.242 & 0.232 & 0.120 & 0.274 \\
    \midrule
    $\lambda_3$ (\textbf{isp}) & 0.270 & 0.200 & 0.248 & 0.257 & 0.279 & 0.088 \\
    \midrule
    $\lambda_4$ (\textbf{obp}) & 0.310 & 0.300 & 0.260 & 0.216 & 0.238 & 0.146 \\
    \bottomrule
    \end{tabular}
    }
    \caption{Statistics for $\mathbb{\lambda_{GS}}$ and $\mathbb{\lambda_{ACS}}$ from \protect\cite{ACS}.}
    \label{tab:lambdas}
\end{table}
In comparison to learning to score cuts directly, this methodology has clear limitations that should be addressed. The performance of this approach is inherently upper bounded by the performance of the metrics being considered in the scoring function. Additionally, the VCG encoding of a given MILP in \cite{ACS} is agnostic to the set of candidate cuts $\mathcal{C}^k$, unlike the tripartite graph from \cite{look_ahead}. It is also well known that MILP solvers suffer from performance variability when changing minute parameters, for instance as simple as the random seed, see \cite{Mip_var}. For this reason, any trained agent obtained with such a method will try to learn optimal cut selection parameters for a specific solver environment with specific parameters being run on specific hardware.

\subsection{Learning \textit{when} to cut}

\citeauthor{local_cuts},~\citeyear{local_cuts}, focus on applying ML to an old and rather Hamletic cut-related question: \emph{To cut or not to cut}? The authors highlight that there is very little understanding on whether a given MIP instance benefits from using only global cuts or also using local cuts, i.e, to Cut (at the root node only) \& Branch or to Branch \& Cut (at all nodes of the tree). They refer to these two alternatives as \emph{"No Local Cuts"} (NLC) and \emph{"Local Cuts"} (LC), respectively, and demonstrate that if access to a perfect decision oracle was possible, a speed-up of 11.36\% is attainable w.r.t. the average solver runtime on a large subset of MIPLIB 2017 instances for FICO Xpress, a commercial MIP solver.

The authors use SL to identify MIP instances that exhibit clear performance gains with one method over the other. Rather than considering this problem as a binary classification task, it is tackled as a regression task predicting the speed-up factor between LC and NLC. The motivation for this is two-fold; first, the ultimate goal is to improve average solver runtime which is a \emph{numerical metric} rather than a \emph{categorical metric}. Second, there are some instances where this decision has negligible impact on solver runtime which complicates creating class labels in the classification approach.

By utilizing feature engineering inspired by the literature on ML for MILP, the authors represent a MILP instance using a 32-dimensional vector that incorporates both static and dynamic features. Static features refer to those that are solver-independent and closely related to the MILP formulation and combinatorial structure. On the other hand, dynamic features are solver-dependent and provide insight into the solver's current behaviour and understanding of a given MILP at different stages in the solution process. The best results were provided by a random forest (RF) which exhibited a speedup of $8.5\%$ on the train set and $3.3\%$ on the test set, prompting further successful experiments which resulted in the implementation of RF as a default heuristic into the new version of FICO Xpress MIP solver.

\section{Theoretical results}
\label{sec:theory}

The nascent line of theoretical work on the sample complexity of learning to tune algorithms for MILP and other NP-Hard problems~\cite{balcan2021much} has also been applied to cut selection. The setting is as follows: There is an unknown probability distribution over MILP instances of which one can access a random subset, the training set. The ``parameters" being learned are $m$ weights which, when applied to the corresponding constraints in the MILP, generate a single cut. \cite{balcan2022structural} and~\cite{tree_complexity} study the number of training instances required to accurately estimate the expected size of the B\&C search tree for a given parameter setting. The expectation here is over the full, unknown distribution over MILP instances. Theorem 5.5 of~\cite{balcan2022structural} shows that the sample complexity is \textit{polynomial} in $n$ and $m$ for Gomory Mixed-Integer cuts, a positive result which is nonetheless not directly useful in practice.

The authors in \cite{ACS} prove that fixed cut selection weights $\lambda$ for Eq.~\eqref{eqn:ScipScore} do not always lead to finite convergence of the cutting plane method and are hence not optimal for all MILP instances. More specifically, consider the convex combination of two metrics (\textbf{isp} and \textbf{obp}), i.e, $\lambda \in \mathcal{R}$ resulting in a finite discretization of $\lambda$ as well an infinite family of MILP instances together with an infinite amount of family-wide valid cuts that can be generated. The use of a pure cutting plane method with a single cut per round can result in the aforementioned instances not terminating to optimality when using $\lambda$ values in the discretization, whereas they terminate for an infinite number of alternative $\lambda$ values.

\section{Conclusion and Future Directions}
\label{sec:conclusion}

Given the rise in the use of ML in combinatorial and integer optimization and the challenges that still exist within, this survey serves as a starting point for future research in integrating data-driven decision-making for cut management in MILP solvers and other related discrete optimization settings. We summarized and critically examined the state-of-the-art approaches whilst demonstrating why ML is a prime candidate to optimize decisions around cut selection and generation, both through a practical and theoretical lens. 
Given the reliance on many default parameters that do not consider the underlying structure of a given general MILP instance, learning techniques that aim to find instance-specific parameter configurations is an exciting area of future research for both MILP and ML communities.
Many additional future directions remain for the research surrounding Learning to Cut. These can range from revisiting algorithmic configuration for cut-related parameters, using ML to identify new strong scoring metrics, and embedding ML in other cut components such as cut generation or removal. Some challenges arise from the aforementioned discussion on the literature for Learning to Cut:\\


\noindent\textbf{Fair Baselines:} There is a lack of fair baselines for ML methods that appropriately balance solver viability and the computational expense from training and data collection that may or may not warrant the use of complex methods such as ML. For instance, \cite{ACS} clearly motivates instance-specific weights however given the relatively marginal learning capabilities and even smaller generalization capabilities we believe analysis into methods like algorithmic configuration should be considered.\\


\noindent\textbf{Large-scale parameter exploration:} There is an overall lack of non-commercial and publicly available large-scale instance datasets for the assessment of cut parameter configuration. \cite{local_cuts} is a prime example that not only are there still many decisions in MILP solvers that are not truly understood, but ML serves as a prime candidate to learn optimal instance-specific decision-making in complex algorithmic settings. For example, a recent "challenge" paper \cite{contardo2022cutting} shows small experiments on MIPLIB 2010~\cite{koch2011miplib} that go against the common belief among MILP researchers that conservative cut generation in B\&B is preferred.

\clearpage




\bibliographystyle{named}
\bibliography{ijcai23}

\end{document}